\documentclass[reqno,10pt, centertags,draft]{amsart}
\usepackage{amsmath,amsthm,amscd,amssymb,latexsym,upref}
\usepackage{latexsym}
\usepackage{lscape}
\usepackage{url}

\newcommand{\beq}{\begin{equation}}
\newcommand{\enq}{\end{equation}}

\makeatletter
\def\theequation{\@arabic\c@equation}

\newcommand{\bbN}{{\mathbb{N}}}
\newcommand{\bbR}{{\mathbb{R}}}

\newcommand{\bbZ}{{\mathbb{Z}}}
\newcommand{\bbC}{{\mathbb{C}}}

\newcommand{\bbT}{{\mathbb{T}}}

\newcommand{\cE}{{\mathcal E}}
\newcommand{\cF}{{\mathcal F}}
\newcommand{\cG}{{\mathcal G}}

\newcommand{\cW}{{\mathcal W}}

\newcommand{\lb}{\label}

\newcommand{\bu}{{\mathbf u}}
\newcommand{\bv}{{\mathbf v}}
\newcommand{\bk}{{\mathbf k}}

\newcommand{\bw}{{\mathbf w}}

\newcommand{\bY}{{\mathbf Y}}
\newcommand{\bp}{{\mathbf p}}
\newcommand{\bq}{{\mathbf q}}
\newcommand{\bx}{{\mathbf x}}
\newcommand{\by}{{\mathbf y}}

\newcommand{\bz}{{\mathbf z}}

\newcommand{\spec}{\operatorname{Spec}}

\newcommand{\ran}{\text{\rm{ran}}}

\newcommand{\bi}{\bibitem}

\numberwithin{equation}{section}

\renewcommand{\div}{\operatorname{div}}
\renewcommand{\det}{\operatorname{det}}

\newcommand{\re}{\operatorname{Re}}

\newcommand{\curl}{\operatorname{curl}}

\newcommand{\diag}{\operatorname{diag}}

\theoremstyle{plain}
\newtheorem{theorem}{Theorem}[section]

\newtheorem{lemma}[theorem]{Lemma}
\newtheorem{corollary}[theorem]{Corollary}

\theoremstyle{definition}

\newtheorem{remark}[theorem]{Remark}

\begin{document}
\allowdisplaybreaks

\title[Continued fractions and the Evans function]{Fredholm determinants, continued fractions, Jost and Evans functions for a Jacobi matrix associated with the 2D-Euler equations}

\author[Y. Latushkin]{Yuri Latushkin}
\thanks{Y.L was supported by the NSF grant DMS-2106157, and would like to
thank the Courant Institute of Mathematical Sciences at NYU and
especially Prof.\ Lai-Sang Young for their hospitality.}
\address{University of Missouri, Columbia, MO 65211, USA}
\email{latushkiny@missouri.edu}

\author[S. Vasudevan]{Shibi Vasudevan}
\address{Krea University, Sri City, Andhra Pradesh, 517646, India}
\email{shibi.vasudevan@krea.edu.in}
\date{\today}
\dedicatory{To Ilya Spitkovsky on his seventieth birthday}
\keywords{2D Euler equations, unidirectional flows, continued fractions, Jost function, Fredholm determinants, Evans function, unstable eigenvalues}

\begin{abstract}
For a second order difference equation that arises in the study of stability of unidirectional (generalized Kolmogorov) flows for the Euler equations of ideal fluids on the two dimensional torus,  we relate the following five functions of the spectral parameter: the Fredholm determinants of the Birman-Schwinger operator pencils associated with the second order equation and the equivalent system of the first order equations; the Jost function constructed by means of the Jost solutions of the second order equation; the Evans function constructed by means of the matrix valued Jost solutions of the first order system, and, finally, to  backward and forward continued fractions associated with the second order difference equation. We prove that all five functions are equal, and that their zeros are the discrete eigenvalues of the second order difference equation. We use this to improve an instability result for a generalization of the Kolmogorov (unidirectional)  flow for the Euler equation on the 2D  torus. 
\end{abstract}

\maketitle

\section{Introduction and the main results} We study the eigenvalue problem \eqref{eulerev1} for a second order asymptotically autonomous difference equation, arising in stability analysis of the generalized Kolmogorov flow of the Euler equation of ideal fluids on 2D torus,
\begin{equation}\label{eulerev1}
z_{n-1}-z_{n+1}=\frac{\lambda}{\rho_n}z_n, \, n\in\bbZ, \,  \lambda\in\bbC, \,\bz:=(z_n)_{n\in\bbZ},
\end{equation}
where $(\rho_n)_{n\in\bbZ}$ is a given sequence satisfying the following conditions:
\begin{equation}\label{condrho}
\rho_0<0,\, \rho_n\in(0,1), \, n\in\bbZ\setminus\{0\}, \text{ and $\rho_n=1+O(1/n^2)$ as $|n|\to\infty$.}
\end{equation}
Sometimes it is convenient to re-write \eqref{eulerev1} as either
\begin{equation}\label{ev3}
\big(S-S^*+\diag_{n\in\bbZ}\{b_nc_n\}\big)\bz=\lambda\bz,
\end{equation}
where $S:(w_n)\mapsto(w_{n-1})$ and $S^*:(w_n)\mapsto(w_{n+1})$ are the shift operators,  or as a first order $(2\times 2)$-system of difference equations,
\begin{equation}\label{evsys1}
y_{n+1}=A_n^{\times}y_n,\, y_n=\begin{bmatrix}z_n\\z_{n-1}\end{bmatrix}\in\bbC^2,\, A_n^{\times}=A+B_nC_n\in\bbC^{2\times 2},\, n\in\bbZ.
\end{equation}
Here and throughout the paper we use the following notations,
\begin{align}\nonumber
A^\times_n&=\begin{bmatrix}-\lambda/\rho_n&1\\1&0\end{bmatrix},\,
A=A(\lambda)=\begin{bmatrix}-\lambda&1\\1&0\end{bmatrix},\, Q_+=\begin{bmatrix}1&0\\0&0\end{bmatrix},\,  Q_-=\begin{bmatrix}0&0\\0&1\end{bmatrix},\\ B_n&=b_nQ_+,\, C_n=c_nQ_+,\, b_n=-\lambda\sqrt{1-\rho_n}/\rho_n,\,
c_n=\sqrt{1-\rho_n}.\label{Anot}
\end{align}
We call $\lambda$ an eigenvalue for \eqref{eulerev1} if the equation has a nonzero solution $\bz$ in the space $\ell^2(\bbZ; \bbC)$ of complex valued sequences.  The essential spectrum of  \eqref{eulerev1} is equal to $\spec(S-S^*)=[-2{\rm i},2{\rm i}]$, and since we are interested in discrete eigenvalues, we assume throughout the paper that
$\lambda\notin[-2{\rm i},2{\rm i}]$. The eigenvalues of $A(\lambda)$ solve the quadratic equation $\mu^2+\lambda\mu-1=0$, and so the assumption is equivalent to the eigenvalues of $A(\lambda)$ being  off the unit circle.

One can associate with \eqref{eulerev1} five functions of the spectral parameter $\lambda$ which we will now define. We begin with \eqref{ev3} and introduce the following Birman-Schwinger-type pencil 
of operators acting in $\ell^2(\bbZ;\bbC)$,
\begin{equation}\label{defK}
K_\lambda=-\diag_{n\in\bbZ}\{c_n\}(S-S^*-\lambda)^{-1}\diag_{n\in\bbZ}\{b_n\}.
\end{equation}
 By \eqref{condrho} the operator $K_\lambda$ is of trace class and so we may define our first function, the characteristic determinant $\det(I-K_\lambda)$.

Next, we re-write \eqref{evsys1} as 
$\big(S^*-\diag_{n\in\bbZ}\{A_n^{\times}\}\big)\by=0$ for $\by=(y_n)_{n\in\bbZ}$ where with a slight abuse of notation we continue to denote by $S^*$ 
the shift acting in the space $\ell^2(\bbZ; \bbC^2)$ of vector valued sequences. We  introduce the Birman-Schwinger-type pencil 
of operators acting in $\ell^2(\bbZ; \bbC^2)$,
\begin{equation}\label{defT}
T_\lambda=\diag_{n\in\bbZ}\{C_n\}\big(S^*-\diag_{n\in\bbZ}\{A\}\big)^{-1}\diag_{n\in\bbZ}\{B_n\}.
\end{equation}
By \eqref{condrho} the operator $T_\lambda$ is of trace class and so we may define our second function, the characteristic determinant $\det(I-T_\lambda)$.

To define our third function,  we introduce notation $\mu_\pm=\mu_\pm(\lambda)$ for the eigenvalues of $A=A(\lambda)$ such that 
\begin{equation}\label{muineq}
|\mu_+(\lambda)|<1<|\mu_-(\lambda)|\end{equation}
and $P_\pm=P_\pm(\lambda)$ for the Riesz projection for $A$ with $\spec(A\big|_{\ran P_\pm})=\{\mu_\pm\}$.
We consider the $(2\times 2)$-{\em matrix valued Jost solutions} $\bY^\pm=(Y^\pm_n)_{n\in\bbZ}$, $Y^\pm_n\in\bbC^{2\times 2}$,  which are the (unique) solutions to \eqref{evsys1} satisfying
\begin{equation}\label{mvj}
\big\|(\mu_\pm)^{-n}\big(Y^\pm_n-A^nP_\pm\big)\|_{\bbC^{2\times 2}}\to0 \text{ as $n\to\pm\infty$ and $Y^\pm_0=Y^\pm_0P_\pm$.}
\end{equation}
We stress that $\bY^\pm=\bY^\pm(\lambda)$ from \eqref{mvj}  depends on $\lambda$ and introduce our third function, the {\em Evans function}, by the formula
\begin{equation}\label{defcE}
\cE(\lambda)=\det\big(Y_0^+(\lambda)+Y^-_0(\lambda)\big).
\end{equation}
Next, we recall that solutions $\bz^\pm=(z^\pm_n)_{n\in\bbZ}$ of the second order difference equation \eqref{eulerev1}  are called the {\em Jost solutions} provided 
\begin{equation}\label{zpm-asymp}
(\mu_\pm)^{-n}z^\pm_n-1\to0 \text{ as $n\to\pm\infty$}.
\end{equation}
As for the discrete Schr\"odinger equations, the Jost solutions $\bz^\pm=\bz^\pm(\lambda)$ are unique. Introducing notation $\cW(\bu,\bv)_n=(-1)^{n-1}\big(u_{n-1}v_n-u_nv_{n-1}\big)$ for the Wronskian of any two sequences $\bu=(u_n)_{n\in\bbZ}$ and $\bv=(v_n)_{n\in\bbZ}$, we note that the Wronskian $\cW(\bz^+,\bz^-)_n$ of the solutions to  \eqref{eulerev1}  is $n$-independent, and define our fourth function, the {\em Jost function}, by the formula
\begin{equation}\label{defF}
\cF(\lambda)=\big(\mu_+(\lambda)-\mu_-(\lambda)\big)^{-1}\cW(\bz^+(\lambda),\bz^-(\lambda))_0.
\end{equation}

Our fifth function uses continued fractions associated with 
\eqref{eulerev1}: we introduce notation 
\begin{align}\lb{deffg}
g_+(\lambda)=
\cfrac{1}{
\frac{\lambda}{\rho_{1}}+\cfrac{1}{
\frac{\lambda}{\rho_{2}}+\dots}},\,\,
g_-(\lambda)=\cfrac{1}{
\frac{\lambda}{\rho_{-1}}+\cfrac{1}{
\frac{\lambda}{\rho_{-2}}+\dots}}.
\end{align}
The continued fractions converge for $\re(\lambda)>0$ and $|\arg(\lambda)|\le \pi/2-\delta$ for any $\delta\in(0,\pi/2)$ by the classical Van Vleck Theorem \cite[Theorem 4.29]{JT} and we may now define our fifth function of interest by the formula
\begin{equation}\label{defG}
\cG(\lambda)= z_0^+(\lambda)z_0^-(\lambda)\big(\mu_-(\lambda)-\mu_+(\lambda)\big)^{-1}\big(\frac{\lambda}{\rho_0}+g_+(\lambda)+g_-(\lambda)\big),
\end{equation}
where $z^\pm_0(\lambda)$ are the $0$-th entries of the Jost solutions, $\mu_\pm(\lambda)$ are the eigenvalues of $A(\lambda)$, and  $g_\pm(\lambda)$ are defined in \eqref{deffg},
and we assume that $\re(\lambda)>0$.

We are ready to formulate the following simple but useful assertion which is the main result of this paper.

\begin{theorem}\label{thm:main} Assume \eqref{condrho} and  $\lambda\notin[-2{\rm i}, 2 {\rm i}]$. The functions introduced in \eqref{defK}, \eqref{defT}, \eqref{defcE}, \eqref{defF} and \eqref{defG} are  holomorphic in $\lambda$ and equal, 
\begin{equation}\label{bigfive}
\det(I-K_\lambda)=\det(I-T_\lambda)=\cE(\lambda)=\cF(\lambda)=\cG(\lambda),
\end{equation}
where the last equality holds under the additional assumption $\re(\lambda)>0$.
Moreover, $\lambda$ is a discrete eigenvalue of \eqref{eulerev1} if and only if it is a zero of each of the functions in \eqref{bigfive}.
\end{theorem}

\begin{remark}
As the proof of Theorem \ref{thm:main} shows, the part of the theorem unrelated to continued fractions holds under much weaker assumptions than \eqref{condrho}. Indeed, let us consider two {\em arbitrary} sequences $(b_n)_{n\in\bbZ}$ and $(c_n)_{n\in\bbZ}$ (that is, {\em not} necessarily defined as in \eqref{Anot} using a given sequence $(\rho_n)_{n\in\bbZ}$ satisfying \eqref{condrho}). Instead of \eqref{condrho} we assume that $(b_n)_{n\in\bbZ},(c_n)_{n\in\bbZ}\in\ell^2(\bbZ;\bbC)$ and instead of \eqref{eulerev1} we consider a more general difference equation \eqref{ev3}. We then define $B_n$ and $C_n$ as indicated in \eqref{Anot}, and consider the system \eqref{evsys1}. Next, assuming $\lambda\notin[-2 {\rm i}, 2 {\rm i}]$, we define $K_\lambda$, $T_\lambda$, $\cE(\lambda)$ and $\cF(\lambda)$ as above, but do not consider continued fractions from \eqref{deffg}.  In this setting all assertions in Theorem \ref{thm:main} hold for equation \eqref{ev3} instead of \eqref{eulerev1}  except the last equality in \eqref{bigfive}.
\hfill$\Diamond$\end{remark}


We now briefly review the relevant literature. The Birman-Schwinger operator pencils are commonly used tools in scattering theory both abstractly and in the context of continuous and discrete Schr\"odinger equations,
\begin{equation}\label{schreq}
u_{n-1}+u_{n+1}+V_nu_n=\lambda u_n, n\in\bbZ, \,
-\partial_{xx}u+V(x)u(x)=\lambda u(x), x\in\bbR.
\end{equation}
Instead of supplying a long list of references, we refer the reader to \cite{BtEG,CGLNSZ} where one can find further bibliography. The (Fredholm) determinants are often used to detect spectra of operators, the standard references are \cite{GGKr,S05}. The most relevant to the current setup papers are \cite{CL07,GLM07} and \cite{LV18}. 

The Birman-Schwinger type pencils specific  for the first order  systems of difference equations have been studied in \cite{CL07} which is a companion paper to \cite{GLM07} dealing with differential equations. To a large extent these papers are based on the theory of integral operators with semi-separable kernels \cite{GM04,GGK90}. In particular, among other things,  the first equality in \eqref{bigfive} was proved in \cite{GLM07} for the case of Schr\"odinger differential operators. We were not aware of results of this sort for the difference equations of type \eqref{eulerev1}. This equation is of course a particular case of eigenvalue problems for three diagonal (Jacobi) matrices, and we mention important contributions in this topic starting with \cite{JN98} and continued by the same authors. 

The Evans function is a powerful tool to study stability of special solutions of partial differential equations, and we refer to \cite{AGJ,KP,PW1,S} where one can find many other sources. The most relevant papers are again \cite{CL07,GLM07}, where, in particular, the second equality in \eqref{bigfive} has been proved for quite general first order systems of difference and differential equations. Involvement of the Evans function in the study of the 2D-Euler-related difference equations as in \eqref{eulerev1} is relatively new, we are aware of only one
paper \cite{DM24} on this topic. The Jost solutions and the Jost function are classical tools in scattering theory for the Schr\"odinger equations \cite{CS,Ne02}. The equality of $\det(I-K_\lambda)$ and $\cF$ in \eqref{bigfive} for the Schr\"odinger case is a classical result going back to \cite{JP51}, see also \cite{Ne72,S05}, that has been extended in many different directions, see \cite{GLZ08,GMZ07,GN15} and the literature therein; however, it is quite possible that the use of the Jost solutions and Jost function as well as the equality of $\cF$ and $\cE$ and $\det(I-K_\lambda)$ in the context of \eqref{eulerev1} appear to be new.

The continued fractions are being used in the context of this equation since probably 
\cite{MS} and we refer to \cite{DLMVW20,FH98,JT} for further references. In particular,  the proof of the last equality in \eqref{bigfive} is based on \cite{DLMVW20}.

We now address in some details the connections of Theorem \ref{thm:main} to hydrodynamic instabilities. 
Indeed, although \eqref{eulerev1} looks very particular, the equations of this type are important as they arise in stability issues for the steady states solutions of many equations of hydrodynamics. Leaving the ``many'' part of the last sentence to a forthcoming paper, we now explain how \eqref{eulerev1} appears in just one particularly important case, that is, in the study of stability of the so-called unidirectional, or generalized Kolmogorov, flow for the Euler equations of ideal fluid on the two dimensional torus $\bbT^2$,
\begin{equation}\label{e}
\partial_{t} \Omega + U \cdot \nabla \Omega = 0,
\, \div U=0, \, \Omega=\curl U, \,  \bx=(x_1,x_2)\in\bbT^2,
\end{equation}
where the two-dimensional vector $U=U(\bx)$ is the velocity and the scalar $\Omega$ is the vorticity of the fluid. Using the Fourier series $\Omega(\bx)=\sum_{\bk\in\bbZ^2\setminus\{0\}}\omega_\bk {\rm e}^{{\rm i} \bk\cdot\bx}$ for vorticity, we rewrite \eqref{e} as
\begin{equation}\label{DEEalpha}
\frac{d\omega_{\mathbf k}}{dt}=\sum_{\mathbf q\in\mathbb Z^2\setminus\{0\}}\beta(\mathbf k-\mathbf q,\mathbf q)\omega_{\mathbf k-\mathbf q}\omega_\mathbf q,\,\,\mathbf k\in\mathbb Z^2\setminus\{0\},
\end{equation}
where the coefficients $\beta(\mathbf p,\mathbf q)$ for $\mathbf p,\mathbf q\in\mathbb Z^2$ are defined as 
\begin{equation}\label{dfnalpha}
\beta(\mathbf p,\mathbf q)=\frac{1}{2}\bigg(\|\mathbf q\|^{-2}-\|\mathbf p\|^{-2}\bigg)(\mathbf p\wedge\mathbf q)\,
\end{equation}
for 
$\mathbf p\neq0,\mathbf q\neq0$, and $\beta(\mathbf p,\mathbf q)=0$ otherwise. Here 
\begin{equation}\label{dfnwedge}
\mathbf p\wedge\mathbf q=\det\left[\begin{smallmatrix}p_1&q_1\\ p_2&q_2\end{smallmatrix}\right] \mbox{ for } \mathbf p=(p_1,p_2) \mbox{ and } \mathbf q=(q_1,q_2).
\end{equation}

The unidirectional (or generalized Kolmogorov) flow is defined as follows: we
fix $\mathbf p \in \mathbb Z^{2} \backslash \{0\}$ and $\alpha\in\bbR$ and consider the steady state solution to the Euler equations on $\mathbb T^{2}$ of the form
\begin{equation}\label{uni}
\Omega^{0}(\mathbf x) = \alpha e^{i \mathbf p \cdot \mathbf x}/2+ \alpha e^{-i \mathbf p \cdot \mathbf x}/2 = \alpha \cos (\mathbf p \cdot \mathbf x).
\end{equation}
In particular, $\bp=(m,0)\in\bbZ^2$ for $m\in\bbN$ corresponds to the classical Kolmogorov flow.
It is thus a classical problem to study (linear) stability of the flow given by \eqref{uni}. To this end, we linearize \eqref{DEEalpha} about the unidirectional flow and consider in $\ell^2(\mathbb Z^2;\bbC)$ the following operator,
\begin{align}
L&: (\omega_\mathbf k)_{\mathbf k\in\mathbb Z^2}\mapsto
\big(\alpha\beta(\mathbf p,\mathbf k-\mathbf p)\omega_{\mathbf k-\mathbf p}-
\alpha\beta(\mathbf p,\mathbf k+\mathbf p) \omega_{\mathbf k+\mathbf p}\big)_{\mathbf k\in\mathbb Z^2}.\label{dfnLB}
\end{align}
The flow \eqref{uni} is called (linearly) unstable if the operator $L$ has
 nonimaginary spectrum. Notice that there are general results for the $2D$-Euler equation  that  assert {\em nonlinear} instability in appropriate function spaces in  the presence of unstable eigenvalues, cf., e.g., \cite{ZL}. See \cite{BGS02, FSV97, HMRW85, LV19, ZL04} for further results on nonlinear stability and instability for the 2D Euler equations and related models.
 
To study the spectrum of $L$ we decompose this operator 
 into a sum of operators $L_{\mathbf q}$, $\mathbf q\in\mathbb Z^2$, acting in the space $\ell^2(\mathbb Z;\bbC)$ by ``slicing'' the grid $\bbZ^2$ along lines parallel to $\bp$. Indeed, for each $\mathbf q\in\mathbb Z^2$ the operator $L$ from \eqref{dfnLB} leaves invariant the subspace of sequences $(\omega_\bk)_{\bk\in\bbZ^2}$ supported on the set $B_\bq=\{\bk=\mathbf q+n\bp: n\in\bbZ\}\subset\bbZ^2$.  The restriction of $L$ to the subspace gives rise to the operator $L_\bq$ defined in $\ell^2(\mathbb Z;\bbC)$ as follows,
\begin{equation}
L_{\mathbf q}: (w_n)\mapsto
\big(\alpha\beta(\mathbf p,\mathbf q+(n-1)\mathbf p) w_{n-1}-
\alpha\beta(\mathbf p,\mathbf q+(n+1)\mathbf p) w_{n+1}\big),\label{dfnLB1}
\end{equation}
where $n\in\mathbb Z$ and for $\mathbf k=\mathbf q+n\mathbf p$ from \eqref{dfnLB} we denote $w_n=\omega_{\mathbf q+n\mathbf p}$. In other words,
 $L_\bq=\alpha(S-S^*)\diag\{\beta(\bp, \bq+n\bp)\}$.
By \eqref{dfnalpha}, if $\mathbf q\parallel\bp$ then $L_{\mathbf q}=0$ and thus we assume throughout that $\mathbf q\not{\parallel}\mathbf p$. Clearly, $\spec(L)=\cup_\bq\spec(L_\bq)$.

To simplify exposition,  in this paper we consider only the case when 
$\bp$ and $\bq$ are such that 
\begin{equation}\label{condpq}
\|\bq\|<\|\bp\| \text{ and $\|\bq+n\bp\|>\|\bp\|$ for all $n\in\bbZ\setminus\{0\}$.}\end{equation} Clearly, a $\bq$ with this property does {\em not} exist for all $\bp$ (say, for $\bp=(1,0)$, corresponding to the Kolmogorov flow known to be stable, or for $\bp=(1,2)$). Others than \eqref{condpq} possibilities for $\bq$ and $\bp$ will be treated elsewhere; our current
assumption corresponds to the case $I_0$ described in \cite{DLMVW20}, and we refer to this paper for a discussion regarding other possible cases.
Moreover, since  $L_\bq$ contains a scalar multiple $\alpha\in\bbR$, with no loss of generality we may rescale this operator or, equivalently, may assume the normalization condition
\begin{equation}\label{normcond}
\alpha(\bq\wedge\bp)\|\bp\|^{-2}/2=1. 
\end{equation}
 
 With these assumptions and conventions, we introduce the sequence 
$\rho_n=1-\|\bp\|^2\|\bq+n\bp\|^{-2}$, $n\in\bbZ$, such that the operator $L_\bq$ in \eqref{dfnLB1} reads  $L_\bq=(S-S^*)\diag\{\rho_n\}$, see \eqref{dfnalpha}. Then condition 
 \eqref{condrho} holds due to \eqref{condpq}. Setting $z_n=\rho_nw_n$, we rewrite the eigenvalue problem $L_\bq\bw=\lambda\bw$ for the operator $L_\bq$, that is, $\rho_{n-1}w_{n-1}-\rho_{n+1}w_{n+1}=\lambda w_n$ as equation \eqref{eulerev1}. Thus, Theorem \ref{thm:main} applies, and we
 are ready to present its consequences.
  \begin{corollary}\label{cor:main}
 Assume that a given $\bp\in\bbZ^2$ is such that there exists a $\bq\in\bbZ^2$ satisfying
 \eqref{condpq}. Then the eigenvalues of $L_\bq$ with positive real parts are in one-to-one  correspondence with  zeros of each of the five functions in \eqref{bigfive}. As a result, the unidirectional flow \eqref{uni} is linearly unstable because the operator $L_\bq$ from \eqref{dfnLB1} and thus $L$ from \eqref{dfnLB} have a positive eigenvalue.
 \end{corollary}
 The instability of the unidirectional flow in the current setting has been established in \cite[Theorem 2.9]{DLMVW20}. Nevertheless, the first assertion in Corollary \ref{cor:main}  is an improvement of the part of \cite[Theorem 2.9]{DLMVW20} where the correspondence was established between merely the positive roots of the function $\lambda/\rho_0+g_+(\lambda)+g_-(\lambda)$ and the positive eigenvalues of $L_\bq$ and under the additional assumption that the respective eigensequences satisfy some special property \cite[Property 2.8]{DLMVW20}. By applying Theorem \ref{thm:main} in the proof of Corollary \ref{cor:main} we were able to show that this assumption is redundant.
 
 We conclude this section with references on the literature on stability of unidirectional flows. This topic is quite classical and well-studied, and we refer to \cite{BW,BFY99,FZ98,FH98,FSV97,FVY00,MS}. The setup as in \eqref{DEEalpha} was also used in many papers \cite{DLMVW20, WDM1, WDM2, LLS, L, SV21}, but the closest to the current work is \cite{DLMVW20}. We mention also \cite{LV18} regarding connections to the Birman-Schwinger operators. Finally, connections between the Evans function and the linearization of the 2D Euler equation has been studied in a recent important paper \cite{DM24}.

\section{Proofs of the main results}

We begin with several comments regarding the objects introduced above.

\begin{remark}\label{rem:zero}
The appearance of the operator $K_\lambda$ from \eqref{defK} is quite natural: Indeed, we factor out to the left the resolvent of $S-S^*$ in \eqref{ev3}, and use that the nonzero elements of the spectra of the operators $I+BC$ and $I+CB$ are the same to conclude that a nontrivial solution $\bz\in\ell^2(\bbZ;\bbC)$ of \eqref{ev3} exists if and only if the operator $I-K_\lambda$ has a nonzero null-space. 
By the assumption $\lambda\notin[-2{\rm i}, 2 {\rm i}]$, cf.\ \eqref{muineq}, the constant coefficient operator $S^*-\diag_{n\in\bbZ}\{A\}$ (from now on we will write this operator as $S^*-A$) is invertible in $\ell^2(\bbZ;\bbC^2)$. Writing 
\begin{align*}
S^*-\diag_{n\in\bbZ}\{A_n^{\times}\}=\big(S^*-A\big)
\big(I-\big(S^*-A\big)^{-1}\diag_{n\in\bbZ}\{B_n\}\diag_{n\in\bbZ}\{C_n\}\big)\end{align*}
and passing from the operator $I+BC$ to $I+CB$ as before we conclude that \eqref{evsys1} has a nontrivial solution $\by\in\ell^2(\bbZ;\bbC^2)$ if and only if the null-space of the operator $I-T_\lambda$ is nonzero. To sum things up, we see that $\lambda$ is a discrete eigenvalue of \eqref{eulerev1} if and only of $\det(I-T_\lambda)=\det(I-K_\lambda)=0$.\hfill $\Diamond$ \end{remark}

We recall from \cite{CL07} and \cite{GLM07} that the unique matrix-valued Jost solutions $\bY^\pm$ are obtained as solutions to the following $(2\times 2)$-matrix Volterra equations,
\begin{equation}\label{jost1}\begin{split}
Y^+_n-A^nP_+&=-\sum_{k=n}^{+\infty}A^{n-k-1}B_kC_kY_k^+,\\
Y^+_n-A^nP_-&=\sum^{n-1}_{k=-\infty}A^{n-k-1}B_kC_kY_k^-.
\end{split}
\end{equation} 
The existence and uniqueness of $\bY^\pm$ follows by passing to the new unknowns $\widehat{Y}_n^\pm=(\mu_\pm)^{-n}Y_n^\pm$ and using a uniform contraction acting in the spaces of bounded sequences; we refer the reader to the proofs of Theorem 3.2 and 6.5  and Corollary 6.6 of \cite{GLM07} and \cite[Lemma 4.2]{CL07} as well as to Remark \ref{holF} below.

The Jost solutions $\bz^\pm=(z_n^\pm)_{n\in\bbZ}$ are obtained as solutions to  the following scalar Volterra equations,
\begin{equation}\label{jostz2-new}\begin{split}
z_n^+-(\mu_+)^n&=-(\mu_+-\mu_-)^{-1}\sum_{k=n}^{+\infty}b_kc_k\big((\mu_+)^{n-k}
-(\mu_-)^{n-k}\big)z^+_k,\\
z_n^--(\mu_-)^n&=(\mu_+-\mu_-)^{-1}\sum_{k=-\infty}^{n-1}b_kc_k\big((\mu_+)^{n-k}
-(\mu_-)^{n-k}\big)z^-_k.
\end{split}
\end{equation}
A computation using $(\mu_\pm)^2+\lambda\mu_\pm-1=0$ shows that $z_n^\pm$ from \eqref{jostz2-new} satisfy \eqref{eulerev1}, see \cite{LG21} for a similar compuation.
In what follows we will need solutions $z_n^\pm$ to \eqref{jostz2-new} only for $\pm n\ge0$ but since solutions of \eqref{jostz2-new} solve \eqref{eulerev1}, they can be propagated to be defined for all $n\in\bbZ$.
We record the following useful properties of the Jost solutions.

\begin{remark}\label{holF} The Jost solutions $\bz^\pm(\lambda)$ are unique and holomorphic in $\lambda$ by the following standard argument, cf.\ \cite{CL07,GLM07} and the references therein. Letting $h(n)=(\mu_+-\mu_-)^{-1}((\mu_+/\mu_-)^n-1)$ and passing in \eqref{jostz2-new} to the new unknowns $\widehat{z}_n^\pm=(\mu_\pm)^{-n}z_n^\pm$, one obtains equations
\begin{equation}\label{tvo}
\widehat{z}_n^+-1=\sum_{k=n}^{+\infty}b_kc_kh(k-n)\widehat{z}_k^+\text{ and }\,
\widehat{z}_n^--1=\sum_{k=-\infty}^{n-1}b_kc_kh(n-k)\widehat{z}_k^-.
\end{equation}
Since $(b_kc_k)_{k\in\bbZ}\in\ell^1(\bbZ;\bbC)$ by \eqref{condrho}, and using \eqref{muineq}, one checks that the right hand sides of the equations in \eqref{tvo} define uniform contractions in the spaces $\ell^\infty([N,+\infty)\cap\bbZ; \bbC)$ and $\ell^\infty((-\infty, -N]\cap\bbZ; \bbC)$, respectively,  provided $N\in\bbN$ is chosen sufficiently large.
Since $b_n$ are linear in $\lambda$ by \eqref{Anot}, one can choose $N$ uniformly in $\lambda$ on compact sets. Solving \eqref{tvo}, one obtains for $\pm n\ge N$ holomorphic in $\lambda$ solutions $z^\pm_n=z^\pm_n(\lambda)$ of \eqref{jostz2-new} satisfying \eqref{zpm-asymp}. Then $z^\pm_n$ solve \eqref{eulerev1}, and as such can be propagated for $\pm n<N$.
\hfill $\Diamond$ \end{remark}

\begin{lemma}\label{znpm0} If $\lambda>0$ then the solutions $\bz^\pm=(z_n^\pm)$ satisfy $z_n^\pm\neq0$ for $\pm n\ge0$.
\end{lemma}
\begin{proof}
Arguing by contradiction, we suppose that, say, $z_0^+=0$ (the cases when $z_n^+=0$ for some $n>0$ or $z_n^-=0$ for some $n\le0$ is considered analogously). Then \eqref{eulerev1} with $n=1$ yields $z_2^+=-\lambda z_1^+/\rho_1$. Since both $\lambda$ and $\rho_1$ are positive, $z_1^+$ and $z_2^+$ have opposite signs. Clearly, $z_1^+$ cannot be zero as otherwise $z_n^+=0$ for all $n\ge0$ contradicting \eqref{zpm-asymp}. By induction, \eqref{eulerev1} and $\rho_n>0$ yield that $z_n^+$ and $z_{n+1}^+$ have opposite signs and are not zero for all $n\ge1$. We now introduce a positive sequence $w_n=-z_{n-1}^+/z^+_n$, $n\ge2$, satisfying the equation $w_{n+1}=(w_n+\lambda/\rho_n)^{-1}$ for all $n\ge2$ due to \eqref{eulerev1}. The sequence $(w_n)_{n\ge2}$ is bounded and separated from zero: Since $\rho_n\to1$ as $n\to\infty$ by \eqref{condrho},  for large $n$ we have $(w_{n+1})^{-1}=w_n+\lambda/\rho_n>\lambda/2$ and 
so $w_{n+1}<2/\lambda$ yielding $(w_{n+1})^{-1}=w_n+\lambda/\rho_n<2/\lambda+\lambda/2$. Letting $\overline{w}=\limsup_{n\to\infty} w_n$ and $\underline{w}=\liminf_{n\to\infty} w_n$ gives
\[\overline{w}=\limsup_{n\to\infty} w_{n+1}=\limsup_{n\to\infty}(w_n+\lambda/\rho_n)^{-1}=(\underline{w}+\lambda)^{-1};\]
likewise, $\underline{w}=(\overline{w}+\lambda)^{-1}$
which in turn implies $\overline{w}=\underline{w}=\mu_+(\lambda)\in(0,1)$, see \eqref{defmu}. Rewriting \eqref{zpm-asymp} as $z_n^+=(\mu_+)^n(1+o(1))$, $n\to\infty$,  yields
\begin{equation*}
\mu_+(\lambda)=\lim_{n\to\infty} w_n=-\lim_{n\to\infty}\frac{z^+_{n-1}}{z^+_{n}}=-\lim_{n\to\infty}\frac{(\mu_+)^{n-1}(1+o(1))}{(\mu_+)^n(1+o(1))}=-\frac1{\mu_+(\lambda)},
\end{equation*}
a contradiction completing the proof.
\end{proof}
\begin{remark}
As in standard scattering theory, one can show that $\lambda$ is an eigenvalue of the operator $L_{\mathbf q}$ if and only if the Jost solutions $\mathbf z^{\pm}$ are proportional to the eigensequence $\mathbf z$. \hfill$\Diamond$
\end{remark}
We will now discuss how the matrix valued Jost solutions $\bY^\pm$ for \eqref{evsys1}
are related to the Jost solutions $\bz^\pm$ for \eqref{eulerev1}. For this and other purposes it is convenient to diagonalize $A$. We introduce matrices 
\begin{align}
W&=\begin{bmatrix}\mu_+&\mu_-\\1&1\end{bmatrix},\,
W^{-1}=\frac1{\mu_+-\mu_-}\begin{bmatrix}1&-\mu_-\\-1&\mu_+\end{bmatrix},\nonumber\\
\widetilde{A}&=\diag\{\mu_+,\mu_-\},\, M=\frac1{\mu_+-\mu_-}\begin{bmatrix}\mu_+&\mu_-\\-\mu_+&-\mu_-\end{bmatrix}\label{winv} 
\end{align}
so that one has, cf.\ \eqref{Anot},
\begin{equation}\label{AA}
W^{-1}AW=\widetilde{A},\, W^{-1}P_\pm W=Q_\pm, \,
W^{-1}Q_+W=M.
\end{equation}
We now introduce the new unknowns $\widetilde{Y}_n^\pm=W^{-1}Y_n^\pm W$ for which  the equations in \eqref{jost1} become
\begin{equation}\label{jost1t}\begin{split}
\widetilde{Y}_n^+-\widetilde{A}^nQ_+&=-\sum_{k=n}^{+\infty}b_kc_k\widetilde{A}^{n-k-1}M\widetilde{Y}_k^+,\\
\widetilde{Y}_n^--\widetilde{A}^nQ_-&=\sum_{k=-\infty}^{n-1}b_kc_k\widetilde{A}^{n-k-1}M\widetilde{Y}_k^-.
\end{split}
\end{equation}
The second condition in \eqref{mvj} yields
$\widetilde{Y}_n^\pm=\widetilde{Y}_n^\pm Q_\pm$. The last equality means that one of the  columns of $\widetilde{Y}_n^\pm$ is zero. We denote the other column by $\widetilde{y}_n^\pm\in\bbC^{2\times 1}$ so that 
$\widetilde{Y}_n^+=\begin{bmatrix}\widetilde{y}_n^+&0_{2\times 1}\end{bmatrix}$
and $\widetilde{Y}_n^-=\begin{bmatrix}0_{2\times 1}&\widetilde{y}_n^-\end{bmatrix}$. By \eqref{jost1t}, the sequences  $(\widetilde{y}_n^\pm)$ solve the equations
\begin{equation}\label{tildet2}\begin{split}
\widetilde{y}_n^+-\begin{bmatrix}(\mu_+)^n&0\end{bmatrix}^\top&=-\sum_{k=n}^{+\infty}b_kc_k\widetilde{A}^{n-k-1}M\widetilde{y}_k^+,\\
\widetilde{y}_n^--\begin{bmatrix}0&(\mu_-)^n\end{bmatrix}^\top&=\sum_{k=-\infty}^{n-1}b_kc_k\widetilde{A}^{n-k-1}M\widetilde{y}_k^-,
\end{split}\end{equation}
where $\top$ denotes transposition. We now introduce 
\begin{equation}\label{defy}
y_n^\pm=(\mu_\pm)^{-1}W\widetilde{y}_n^\pm.\end{equation}  Since $W\widetilde{Y}_n^\pm=Y_n^\pm W$, the vector  $y_n^\pm$ is a linear combination of the columns of the matrix $Y_n^\pm$ and therefore the sequence $(y_n^\pm)$ is a solution to the equation \eqref{evsys1} which, in turn, is equivalent to \eqref{eulerev1}. Thus, there is a solution $\bz=(z_n^\pm)$ of  \eqref{eulerev1} such that $y_n^\pm=\begin{bmatrix} z_n^\pm & z_{n-1}^\pm\end{bmatrix}^\top$. We now multiply \eqref{tildet2} from the left by $(\mu_\pm)^{-1}W$, use \eqref{defy}, and directly compute the product \[W\widetilde{A}^{n-k-1}MW^{-1}=\frac1{\mu_+-\mu_-}\begin{bmatrix}\mu_+^{n-k}-\mu_-^{n-k}&0\\\mu_+^{n-k-1}-\mu_-^{n-k-1}&0\end{bmatrix}\]
to obtain for $z_n^\pm$  the desired equation \eqref{jostz2-new}.

We are ready to present the proof of Theorem \ref{thm:main}.

\begin{proof} The four equalities in \eqref{bigfive} are proved as follows.

 {\bf 1.} We prove $\det(I-K_\lambda)=\det(I-T_\lambda)$.\, We claim that $K_\lambda$ from \eqref{defK} is the $(1,1)$-block of the operator $T_\lambda$ from \eqref{defT} in the decomposition \[\ell^2(\bbZ;\bbC^2)=\ran(\diag\{Q_+\})\oplus\ran(\diag\{Q_-\})\] for the projections $Q_\pm$ from \eqref{Anot}, that is, that
\begin{equation*}
\diag\{C_n\}(S^*-A)^{-1}\diag\{B_n\}=\begin{bmatrix}-\diag\{c_n\}(S-S^*-\lambda)^{-1}\diag\{b_n\}&0\\0&0\end{bmatrix}.
\end{equation*}
Clearly, this implies the required equality. To begin the proof of the claim, for any $\bz=(z_n)_{n\in\bbZ}\in\ell^2(\bbR;\bbC)$ and $\bu=(u_n)_{n\in\bbZ}\in\ell^2(\bbR;\bbC^2)$   one directly checks the formulas
\begin{align*}
\big((S-S^*-\lambda)^{-1}\bz\big)_n=-\frac{1}{\mu_+-\mu_-}\big(\sum_{k=n}^{+\infty}
(\mu_-)^{n-k}z_k+\sum_{k=-\infty}^{n-1}(\mu_+)^{n-k}z_k\big),\\
\big((S^*-A)^{-1}\bu\big)_n=-\sum_{k=n}^{+\infty}A^{n-k-1}P_-u_k+\sum_{k=-\infty}^{n-1}A^{n-k-1}P_+u_k.
\end{align*}
Using \eqref{AA}, and denoting by $z_n$ the first component of $u_n\in\bbC^2$,
we calculate
\begin{align*}
\big(&\diag\{Q_+\}(S^*-A)^{-1}\diag\{Q_+\}\bu\big)_n=
-\sum_{k=n}^{+\infty}Q_+A^{n-k-1}P_-Q_+u_k\\&\hskip4cm+\sum_{k=-\infty}^{n-1}Q_+A^{n-k-1}P_+Q_+u_k\\&
=-\sum_{k=n}^{+\infty}Q_+W\widetilde{A}^{n-k-1}(I-Q_+)W^{-1}Q_+u_k\\&\hskip4cm
+\sum_{k=-\infty}^{n-1}Q_+W\widetilde{A}^{n-k-1}Q_+W^{-1}Q_+u_k\\&
=-\sum_{k=n}^{+\infty}Q_+W(I-Q_+)W^{-1}Q_+(\mu_-)^{n-k-1}u_k\\&\hskip4cm
+\sum_{k=-\infty}^{n-1}Q_+WQ_+W^{-1}Q_+(\mu_+)^{n-k-1}u_k\\&
=\frac{1}{\mu_+-\mu_-}\begin{bmatrix}\sum_{k=n}^{+\infty}
(\mu_-)^{n-k}z_k+\sum_{k=-\infty}^{n-1}(\mu_+)^{n-k}z_k&0\\0&0\end{bmatrix},
\end{align*}
and now the claim follows by multiplying the last equality by $\diag\{C_n\}$ and $\diag\{B_n\}$, cf.\  \eqref{Anot}.

{\bf 2.} We prove $\det(I-T_\lambda)=\cE(\lambda)$.\, This, in fact,  is proved in \cite[Proposition 5.2]{CL07}  for rather significantly more general than \eqref{evsys1} system of difference equations.

 {\bf 3.} We prove  $\cE(\lambda)=\cF(\lambda)$.\,  Using $\widetilde{Y}_0^\pm=W^{-1}Y_0^\pm W$, formula \eqref{defy}, 
and the relation $y_0^\pm=\begin{bmatrix} z_0^\pm & z_{-1}^\pm\end{bmatrix}^\top$, one infers
\begin{align*}
W^{-1}\big(Y_0^++Y_0^-\big)W&=
\widetilde{Y}_0^++\widetilde{Y}_0^-
=\begin{bmatrix}\widetilde{y}_0^+&\widetilde{y}_0^-\end{bmatrix}
=\begin{bmatrix}\mu_+W^{-1}y_0^+&\mu_-W^{-1}y_0^-\end{bmatrix}\\&
=W^{-1}\begin{bmatrix}z_0^+&z_0^-\\z_{-1}^+&z_{-1}^-\end{bmatrix}\diag\{\mu_+,\mu_-\}.
\end{align*}
Computing determinants, \eqref{winv}, \eqref{defcE} and \eqref{defF}, and using the fact that $\mu_{+}\mu_{-}=-1$, yield the required equality.

 {\bf 4.} We prove $\cF(\lambda)=\cG(\lambda)$.\, Here, we assume that $\re(\lambda)>0$ and so the eigenvalues $\mu_\pm(\lambda)$ in \eqref{muineq} are given by   
\begin{equation}\label{defmu}
\mu_\pm=\mu_\pm(\lambda)=-\frac{\lambda}{2}\pm\sqrt{\big(\frac{\lambda}{2}\big)^2+1}.
\end{equation}
The function $\cF$ is holomorphic in $\lambda$ by Remark \ref{holF}. The convergent continued fractions in \eqref{deffg} are holomorphic in $\lambda$ by the classical Stieltjes-Vitali Theorem \cite[Theorem 4.30]{JT}. Thus, it is enough to show the desired equality only for $\lambda>0$ which we assume from now on.
 We recall that $z_n^\pm\neq0$ for $\pm n\ge0$ by Lemma \ref{znpm0}, introduce $v^\pm_n=z^\pm_{n-1}/z^\pm_n$ and re-write equation 
 \eqref{eulerev1} for $\bz^\pm$ as follows,
 \begin{equation}\label{ueq}
 v^+_n=\frac{1}{\rho_n}+\frac{1}{v^+_{n+1}}\, \text{ for $n\ge0$ and 
 $v^-_{n+1}=-\frac{1}{\frac{\lambda}{\rho_n}-v^-_n}$ for $n\le-1$.}
 \end{equation}
 The first formula can be iterated forward and the second backward producing truncated continued fractions as in \eqref{deffg}.
 Because the continued fractions converge and $\lambda>0$, an argument from \cite[pp.2063]{DLMVW20} involving monotonicity of the sequences formed by the odd and even truncated continued fractions yields 
 \begin{equation}\label{u0}
 v_0^+=\frac{\lambda}{\rho_0}+g_+(\lambda) \text{ and  $v_0^-=-g_-(\lambda)$}.
 \end{equation}
 Therefore,
 \[\frac{\lambda}{\rho_0}+g_+(\lambda)+g_-(\lambda)=v^+_0-v^-_0
 =\frac{z^+_{-1}}{z^+_0}-\frac{z^-_{-1}}{z^-_0}=-\frac{\cW(\bz^+,\bz^-)_0}{z^+_0 z^-_0},\]
 and the desired equality of $\cF$ and $\cG$ follows by \eqref{defF} and \eqref{defG}.
 
This concludes the proof of \eqref{bigfive}. The last assertion in the theorem follows by Remark \ref{rem:zero}.
\end{proof}

We now present the proof of Corollary \ref{cor:main}.

\begin{proof} We claim that the operator $L$ in \eqref{dfnLB} has a positive eigenvalue.
 Since $\spec(L_\bq)\subset\spec(L)$ for a $\bq$ satisfying \eqref{condpq}, it suffices to show that the operator $L_\bq$ from \eqref{dfnLB1} has a positive eigenvalue. Multiplying this operator by a positive factor, we may assume that the normalization condition \eqref{normcond} does hold, and so we need to show that \eqref{eulerev1} has a positive eigenvalue. Now Theorem \ref{thm:main} applies yielding the first assertion in the corollary. To finish the proof it suffices  to show that there is a positive root of the function $\cG$ defined in \eqref{defG}, equivalently, that  for some $\lambda>0$ and $g_\pm$  defined in \eqref{deffg} one has  $-\lambda/\rho_0=g_+(\lambda)+g_-(\lambda)$. We recall that $\rho_0<0$ and $g_\pm(\lambda)>0$ for $\lambda>0$ because condition \eqref{condrho} holds. The proof is completed by using the relations
 $\lim_{\lambda\to0^+}g_\pm(\lambda)=1$ and $\lim_{\lambda\to+\infty}g_\pm(\lambda)=0$ established in \cite[Lemma 2.10(4)]{DLMVW20}.
\end{proof}

\end{document}